\documentclass[12pt,reqno,a4paper]{amsart}
\usepackage{amsmath,amssymb,amsthm,amsaddr}
\usepackage{enumitem}
\usepackage[margin=2.7cm,top=3.5cm,bottom=3.5cm,footskip=1cm,headsep=1.5cm]{geometry}
\usepackage[utf8]{inputenc}

\usepackage{amsthm}
\usepackage[british]{babel}

\usepackage{caption}
\usepackage{subcaption}
\usepackage{float}

\usepackage{tikz}
\usetikzlibrary{decorations.markings}
\usepackage[europeanresistors]{circuitikz}
  \usetikzlibrary{arrows, automata}
\usepackage{siunitx}
\usetikzlibrary{shapes}
\usetikzlibrary{calc}
\usetikzlibrary{spy,arrows}
\usepackage{pgfplots}
\usepackage{pgfplotstable}
\pgfplotsset{compat=newest}

\usepackage{multirow}

\author[H. Egger, K. Schmidt, and V. Shashkov]{H. Egger and K. Schmidt and V. Shashkov}
\address{Department of Mathematics, TU Darmstadt, Germany}

\title[Convolution quadrature methods for coupled dynamical systems]{Multistep and Runge-Kutta convolution quadrature methods for coupled dynamical systems}

\newtheorem{lemma}{Lemma}[section]
\theoremstyle{definition}
\newtheorem{remark}[lemma]{Remark}

\def\dt{\partial_t}
\def\A{\mathcal{A}}
\def\diag{\text{diag}}
\def\CC{\mathbb{C}}
\def\ddt{\frac{d}{dt}}
\def\curl{\text{curl}}

\begin{document}

\begin{abstract}
We consider the efficient numerical solution of coupled dynamical systems, consisting of a small nonlinear part and a large linear time invariant part, possibly stemming from spatial discretization of an underlying partial differential equation.
The linear subsystem can be eliminated in frequency domain and for the numerical solution of the resulting integro-differential algebraic equations, we propose a a combination of Runge-Kutta or multistep time stepping methods with appropriate convolution quadrature to handle the integral terms. 
The resulting methods are shown to be algebraically equivalent to a Runge-Kutta or multistep solution of the coupled system and thus automatically inherit the corresponding stability and accuracy properties. 
After a computationally expensive pre-processing step, the online simulation can, however, be performed at essentially the same cost as solving only the small nonlinear subsystem.
The proposed method is, therefore, particularly attractive, if repeated simulation of the 
coupled dynamical system is required. 
\end{abstract}

\maketitle

\section{Introduction} \label{sec:introduction}

We consider the efficient numerical approximation of coupled linear-nonlinear dynamical systems
described by systems of differential-algebraic or partial-differential algebraic equations.
Our approach is particularly attractive for problems in which
\begin{itemize}
 \item efficient multiple simulation of the system dynamics is required,
 \item the linear subsystem dominates the dimension of the problem while the small nonlinear part dictates the dynamic behavior, and
 \item the coupling between the subsystems takes place via a limited number of ports. 
\end{itemize}
Such problems arise in a variety of applications, e.g., in multibody dynamical systems,
electromechanical devices, or in models of electric power networks. Our research is motivated 
by field-circuit coupled problems arising in the simulation of electronic circuits with {\em electromagnetic elements}.
In this context, the nonlinear subsystem describes the electric circuit while the linear part 
models the electromagnetic fields described by Maxwell's equations \cite{BartelBaumannsSchoeps11,HiptmairSterz05,SchoepsEtAl10}; details will be given in later sections.

The following two standard approaches for the numerical solution of such coupled problems are widely used in practice:
\begin{itemize}
 \item the simultaneous integration of the coupled linear--nonlinear differential-algebraic system by appropriate time-stepping methods;
 \item the representation of the frequency domain response of the linear subsystem by its transfer function, or approximations thereof, and time discretization of the corresponding reduced nonlinear integro-differential algebraic problem.
\end{itemize}
The first approach is well-understood \cite{BenderskayaDeGersemWeilandClemens04,SchoenmakerEtAl16,SchoepsEtAl10} and can be considered to be reliable. 
If the linear subsystems are of large dimension, e.g., when stemming from discretization of partial differential equations, such a holistic approach is, however, computationally very demanding and not suited for repeated online simulation. 
Model order reduction is thus frequently used to compute low-dimensional approximate state space representations for the linear subsystem in an offline stage which then allow for a fast evaluation of the coupled systems in online simulations \cite{Antoulas05,BennerMehrmannSorensen,Roychowdhury99}. 
Multirate time-integration or waveform relaxation can be utilized to further reduce the computational 
complexity \cite{SchoepsEtAl13}. 

The second approach relies on the fast access to the transfer function of the linear subsystem. 
Approriate rational approximations for the transfer function are therefore constructed \cite{GustavsenSemlyen99,AAA18}, 
which can then be incorporated efficiently in the time-domain simulation of the reduced coupled system by recursive convolution \cite{KapurEtAl96,LinKuh92}.

Let us emphasize that the efficient numerical realization of both approaches mentioned above relies on certain low dimensional approximations of the input-output behavior of the linear subsystem, either in state space or in frequency domain. 
Although much research has been devoted to the development of such approximations \cite{Antoulas05,BennerMehrmannSorensen,AAA18}, a systematic evaluation and control of the approximation errors seems to be difficult.

In this paper, we therefore propose a different strategy that allows to conduct online simulations of the problem with the accuracy and stability of a discretization of the fully coupled system at the cost 
of computing solutions to the reduced nonlinear system after elimination of the linear subproblem. 
Similar as the model reduction approaches outlined above, the method consists of a compute intensive 
offline stage, in which a detailed analysis of the linear subsystem is performed, and an efficient 
online phase in which only the reduced nonlinear subproblem has to be integrated.
The resulting method is therefore particularly well suited for repeated online simulation.

\section{Problem description and outline of the approach}

We consider coupled dynamical systems of the general form
\begin{align}
M(y) \dt y + F(y) &= C^\top z, \label{eq:sys1}\\
E \dt z + A z &= B y.         \label{eq:sys2}
\end{align}
Nonlinear equations with leading order term $M(y) \dt y$ replaced by $\dt Q(y)$, which are
frequently encountered in circuit simulation, could be considered with similar arguments. 
For ease of presentation, we utilize trivial initial conditions 
\begin{align}
y(0)=0 \qquad \text{and} \qquad z(0)=0 \label{eq:sys3} 
\end{align}
and we assume that the system is finite dimensional but possibly large;
in particular, the linear subsystem \eqref{eq:sys2} may arise from space discretization of an underlying partial 
differential equation. 
The matrices $M(y)$ and $E$ in front of the time derivatives are not required to be regular, and the system \eqref{eq:sys1}--\eqref{eq:sys2} therefore constitutes a set of coupled differential-algebraic equations. 
We assume that the (perturbation) index is moderate, such that a stable time discretization by appropriate 
single or multistep methods is possible \cite{BrenanCampbellPetzold96,HairerLubichRoche89}. 

For illustration of our ideas, let us consider the time discretization of \eqref{eq:sys1}--\eqref{eq:sys2} by the implicit Euler method with constant step size $\tau$, leading to a recurrence of the form
\begin{align} 
M(y_n) \frac{y_n-y_{n-1}}{\tau} + F(y_n) &= C^\top z_n, \label{eq:euler1}\\
E \frac{z_n-z_{n-1}}{\tau} + A z_n &= B y_n.           \label{eq:euler2} 
\end{align}
Note that a coupled linear-nonlinear algebraic system has to be solved in every time step which is 
prohibitive for a fast online simulation if the dimension of the system is large. 
As we will show in the following, at least the repeated solution of the linear subproblem can be 
completely avoided during online simulation. 

The linearity of equation \eqref{eq:sys2} and an application of the Laplace transform 
allows to express the response of the linear subsystem in frequency domain as  
\begin{align} \label{eq:Ctzh}
C^\top \widehat z(s) = C^\top (s E + A)^{-1} B \widehat y(s) =: \widehat K(s) y(s).
\end{align}
Here and below, we denote by $\widehat y(s)$, $\widehat z(s)$ the Laplace transforms 
of functions $y(t)$, $z(t)$.
The function $\widehat K(s)$ is called the transfer function of the linear subsystem \eqref{eq:euler2}.
Back transformation to time domain
and insertion of the result into \eqref{eq:sys1} then leads to the integro-differential algebraic equation 
\begin{align} \label{eq:idae}
M(y(t)) \dt y(t) + F(y(t)) = \int_0^t K(t-r) y(r) dr.
\end{align}
Let us recall that the kernel function $K(t)$, i.e., the impulse-response of the linear subsystem \eqref{eq:sys2}, is 
determined implicitly by its Laplace transform $\widehat K(s)$ defined above. 

The discretization of related Volterra-integro-differential equations has been investigated
intensively in the literature; see e.g. \cite{Huang08,Lubich83,Matthys76} and the references given there. 
An application of the implicit Euler method to the differential part of \eqref{eq:idae} and 
an appropriate quadrature rule for the integral part leads to the recursion formula
\begin{align} \label{eq:euler3}
M(y_n) \frac{y_n-y_{n-1}}{\tau}  + F(y_n) &= \tau \sum_{k=0}^n \omega_{n-k} y_k.
\end{align}
We will show that, for appropriate choice of the quadrature weights $(\omega_{n})_{n \ge 0}$, the 
numerical solution $(y_n)_{n \ge 0}$ of method \eqref{eq:euler3} coincides with the $y$ component of 
the solution $(y_n,z_n)_{n \ge 0}$ obtained with method \eqref{eq:euler1}--\eqref{eq:euler2}. 
In this sense, the equivalence of the coupled system \eqref{eq:sys1}--\eqref{eq:sys2} with the reduced problem \eqref{eq:idae} thus remains valid after discretization by appropriate strategies. 
As a consequence, also the stability and convergence properties of method \eqref{eq:euler1}--\eqref{eq:euler2} carry over verbatim to method \eqref{eq:euler3}, if the weights $(\omega_n)_{n \ge 0}$ are chosen appropriately. 
The computation of these weights requires a detailed computational analysis of the linear subproblem 
which can, however, be performed in an offline stage. The online simulation can then be achieved efficiently via the reduced scheme \eqref{eq:euler3}.

The above considerations are not limited to the implicit Euler method, but they can be extended 
to single and multistep methods of higher order accuracy. As shown in the work of Lubich and Ostermann \cite{Lubich88a,Lubich88b,LubichOstermann93}, the appropriate definition of the quadrature 
weights can be derived in the framework of convolution quadrature methods.

\bigskip

The remainder of the paper is organized as follows:
In Section~\ref{sec:convolution}, we present the extension of our approach to Runge-Kutta and multistep methods 
and formally derive the corresponding quadrature weights. 
In Section~\ref{sec:implementation}, we discuss the practical computation of the quadrature weights 
and briefly comment on the efficient implementation of the convolution sums. 
In Section~\ref{sec:fieldcircuit}, we discuss in some detail the derivation of linear--nonlinear coupled models arising in the context of the simulation of electric circuit including electromagnetic elements. 
Two particular model problems and numerical results are presented in Section~\ref{sec:numerics}
in order to illustrate our theoretical results.

\section{A review on convolution quadrature} \label{sec:convolution}

We first discuss the discretization of the coupled system \eqref{eq:sys1}--\eqref{eq:sys2} by 
Runge-Kutta and multistep methods and then derive in detail the appropriate definition of the quadrature 
weights $(\omega_n)_{n \ge 0}$ required for the treatment of the integral terms in \eqref{eq:idae}. 
Our presentation is strongly based on the seminal papers \cite{Lubich88a,Lubich88b,LubichOstermann93} 
on convolution quadrature.

\subsection{Multistep methods}

As a first approach, we consider the discretization by multistep methods. For ease of notation, we restrict the presentation to the BDF methods of Gear, which are known to be particularly well suited for the solution of differential algebraic equations \cite{BrenanCampbellPetzold96,HairerLubichRoche89}. 
The discretization of \eqref{eq:sys1}--\eqref{eq:sys2} by the $m$-step BDF formula with constant time step $\tau$, for instance, leads to the recursion
\begin{align} 
M(y_n) \frac{1}{\tau}\sum_{k=0}^m \alpha_k y_{n-m+k} + F(y_n) &= C^\top z_n \label{eq:bdf1} \\
E \frac{1}{\tau}\sum_{k=0}^m \alpha_k z_{n-m+k} + A z_n &= B y_n. \label{eq:bdf2}
\end{align}
Multiplication of equation \eqref{eq:bdf2} by exponentials $\xi^n$ and summation over $n$ further yields
\begin{align} \label{eq:bdf3}
E \frac{\delta(\xi)}{\tau} z(\xi) + A z(\xi) = B y(\xi),
\end{align}
where $y(\xi) = \sum_{n=0}^\infty y_n \xi^n$ and $z(\xi) = \sum_{n=0}^\infty z_n \xi^n$ denote 
the generating functions of the sequences $(y_n)_{n \ge0}$ and $(z_n)_{n \ge 0}$, and  
the characteristic polynomial
\begin{align} \label{eq:bdf4}
\delta(\xi) = \sum_{k=0}^m \alpha_k \xi^{m-k}
\end{align}
is determined by the coefficients of the multistep method under consideration. 
A simple rearrangement of the terms in formula \eqref{eq:bdf3} and using the definition of the transfer function in \eqref{eq:transferK} leads to 
\begin{align} \label{eq:bdf5}
C^\top z(\xi) 
&= C^\top \left(\frac{\delta(\xi)}{\tau} E + A\right)^{-1} B y(\xi)  
 = \widehat K(\frac{\delta(\xi)}{\tau}) y(\xi),
\end{align}
and a formal expansion of the function $k(\xi)=\widehat K(\frac{\delta(\xi)}{\tau})$ 
into a power series further yields
\begin{align} \label{eq:bdf6}
\widehat K\left(\frac{\delta(\xi)}{\tau}\right) = \sum_{n=0}^\infty \omega_n \xi^n.
\end{align}
Note that the coefficients $(\omega_n)_{n \ge 0}$ in this expansion can be found easily, e.g.,  
by comparing derivatives of the two expressions at $\xi = 0$. 
An efficient alternative computation of the coefficients will be discussed in more detail in the next section.

Inserting the above expression and the power series representations for $z(\xi)$ and $y(\xi)$ 
into \eqref{eq:bdf5}, and then applying the Cauchy product formula to the result leads to
\begin{align} \label{eq:bdf7}
\sum_{n=0}^\infty C^\top z_n \xi^n 
&= (\sum_{n=0}^\infty \omega_n \xi^n) \cdot (\sum_{k=0}^\infty y_k \xi^k) 
 = \sum_{n=0}^\infty ( \sum_{k=0}^n \omega_{n-k} y_k) \xi^n. 
\end{align}
A simple comparison of the coefficients in the two power series on the left and right hand side 
of this identity allows to deduce that 
\begin{align} \label{eq:bdf8}
C^\top z_n = \sum_{k=0}^n \omega_{n-k}  y_k.
\end{align}
Using this expression to replace the right hand side in equation \eqref{eq:bdf1} yields the following.
\begin{lemma} \label{lem:bdf}
Let $(y_n,z_n)_{n \ge 0}$ denote a solution of \eqref{eq:bdf1}--\eqref{eq:bdf2} with $y(0)=0$, $z(0)=0$,
and let the quadrature weights $(\omega_n)_{n \ge 0}$ be determined by \eqref{eq:bdf6} with \eqref{eq:bdf4} and \eqref{eq:Ctzh}. Then 
\begin{align} \label{eq:bdf9}
M(y_n) \frac{1}{\tau} \sum_{k=0}^m \alpha_k y_{n-m+k} + F(y_n) &= \sum_{k=0}^n \omega_{n-k} y_k.
\end{align}
\end{lemma}
\noindent
This shows that the $y$ component of the solution for the coupled system \eqref{eq:sys1}--\eqref{eq:sys2} obtained by the BDF method \eqref{eq:bdf1}--\eqref{eq:bdf2} thus coincides with the solution for the integro-differential equation \eqref{eq:idae} computed by the multistep convolution quadrature method \eqref{eq:bdf9}.
\begin{remark}
The coefficients of the $1$-step BDF formula are given by $\alpha_1=1$, $\alpha_0 = -1$, 
and the characteristic polynomial thus reads $\delta(\xi)=1-\xi$. 
The method \eqref{eq:bdf9} then coincides with method \eqref{eq:euler3} outlined in the introduction.
If the quadrature weights are chosen as in \eqref{eq:bdf6}, then the solution $(y_n)_{n \ge 0}$ 
corresponds to that of method \eqref{eq:euler1}--\eqref{eq:euler2}.
\end{remark}

\subsection{Runge-Kutta methods}

We next discuss the time discretization by the Radau-IIA methods, which are also well-suited for the discretization
of differential algebraic equations \cite{HairerLubichRoche89}. 
The time discretization of the coupled system \eqref{eq:sys1}--\eqref{eq:sys2} by an implicit Runge-Kutta method with $s$ stages leads to numerical schemes of the form
\begin{alignat}{2}
y_{n+1} = y_n + \tau \sum_{j=0}^s \beta_j Y_{nj}',  \qquad
z_{n+1} = z_n + \tau \sum_{j=0}^s \beta_j Z_{nj}',  \label{eq:rk1}
\end{alignat}
with stage values $Y_{ni}$, $Z_{ni}$, $i=1,\ldots,s$ defined by
\begin{alignat}{2}
Y_{ni} = y_n + \tau \sum_{j=0}^s \alpha_{ij} Y'_{nj},  \qquad
Z_{ni} = z_n + \tau \sum_{j=0}^s \alpha_{ij} Z'_{nj},  \label{eq:rk2}
\end{alignat}
and slopes $Y'_{nj}$, $Z'_{nj}$, $j=1,\ldots,s$ determined by  
\begin{alignat}{2}
M(Y_{nj}) Y'_{nj} + F(Y_{nj}) = C^\top Z_{nj},   \label{eq:rk3} \\
         E Z'_{nj} + A Z_{nj} = B Y_{nj}.        \label{eq:rk4} 
\end{alignat}
We write $Z_n=[Z_{n1},\ldots,Z_{ns}]^\top$ in the sequel and use similar short hand notation for the
stage values and slopes. 
Similar as in the previous section, we multiply the three equations for the $z$ component 
by $\xi^n$ and sum over all $n$, which leads to
\begin{align}
\xi^{-1} z(\xi) &= z(\xi) + \tau (\beta^\top \otimes I) Z'(\xi),         \label{eq:rk7} \\
Z(\xi) &= \underline 1 \otimes z(\xi) + \tau (\A \otimes I) Z'(\xi).                \label{eq:rk8} \\
(I \otimes E) Z'(\xi) + (I \otimes A) Z(\xi) &= (I \otimes B) Y(\xi).    \label{eq:rk9} 
\end{align}
Recall that $Y(\xi)=\sum_{n=0}^\infty Y_n \xi^n$ is the generating function for the series $\{Y_n\}_{n \ge 0}$ and  $Z(\xi)$, $Z'(\xi)$ are defined similarly. 
Furthermore, $\A_{ij}=\alpha_{ij}$ and  $\beta_j$ are the coefficients of the Runge-Kutta method,
and $\underline 1$ is the constant one vector.
Elimination of $z(\xi)$ and $Z'(\xi)$ via \eqref{eq:rk7}--\eqref{eq:rk8} and substitution
into \eqref{eq:rk9} allows to express the function $Z(\xi)$ as solution of the $\xi$ dependent equation
\begin{align}  \label{eq:rk10}
\left( \frac{\Delta(\xi)}{\tau} \otimes E + I \otimes \A \right) Z(\xi) = (I \otimes B) Y(\xi)
\end{align}
with symbol $\Delta(\xi)$ characterizing the Runge-Kutta method and defined by 
\begin{align} \label{eq:rk10a}
\Delta(\xi) = \left(\frac{\xi}{1-\xi} \underline 1 \beta^\top + \A\right)^{-1}.
\end{align}
A simple rearrangement of the terms in equation \eqref{eq:rk10} then leads to 
\begin{align}
(I \otimes C^\top) Z(\xi) 
&= (I \otimes C^\top) \left(\frac{\Delta(\xi)}{\tau} \otimes E + I \otimes A\right)^{-1} (I \otimes B) Y(\xi) \nonumber\\
&=: \widehat K\left(\frac{\Delta(\xi)}{\tau}\right) Y(\xi).  \label{eq:rk11}
\end{align}
The evaluation of the transfer function $\widehat K(\frac{\Delta(\xi)}{\tau})$ for the matrix valued arguments will be discussed in more detail at the end of the next section.
Similarly as before, we can formally expand $\widehat K(\frac{\Delta(\xi)}{\tau})$ into a power series 
\begin{align} \label{eq:rk12} 
\widehat K\left(\frac{\Delta(\xi)}{\tau}\right) = \sum_{n=0}^\infty W_n \xi^n
\end{align}
with appropriate weight matrices $(W_n)_{n \ge 0}$ that can again be obtained easily, e.g., 
by comparing derivatives at $\xi=0$. 
Using this expression and the power series representations of $Y(\xi)$ and $Z(\xi)$ in \eqref{eq:rk11}, 
applying the Cauchy product formula, 
and comparing the coefficients of the corresponding sequences then leads to
\begin{align} \label{eq:rk13} 
(I \otimes C^\top) Z_n = \sum_{k=0}^n W_{n-k} Y_k, 
\end{align}
which corresponds almost verbatim to the formula \eqref{eq:bdf7} for multistep methods. 
By decomposition of the vectors $Y_n$, $Z_n$, and the weight matrices $W_n$ into their components 
for the individual stages of the Runge-Kutta method, we finally obtain 
\begin{align} \label{eq:rk14}
C^\top Z_{ni} = \sum_{k=0}^n \Big(\sum_{j=1}^s W_{n-k,ij} Y_{nj} \Big), \qquad i=1,\ldots,s 
\end{align}
which can now be inserted into \eqref{eq:rk3} to obtain the following result.
\begin{lemma} \label{lem:rk}
Let $(y_n,z_n)_{n \ge 0}$, $(Y_n,Z_n)_{n \ge 0}$, $(Y'_n,Z'_n)_{n \ge 0}$ be a solution of
\eqref{eq:rk1}--\eqref{eq:rk4} with $y(0)=0$, $z(0)=0$, and let the weights
$(W_n)_{n \ge 0}$ be defined by \eqref{eq:rk10a}--\eqref{eq:rk12}.
Then 
\begin{align}
y_{n+1} &= y_n + \tau \sum_{j=0}^s \beta_j Y_{nj}',  \label{eq:rk15}
\end{align}
with stages and slopes defined by 
\begin{align}
M(Y_{ni}) Y'_{ni} + F(Y_{ni}) = \sum_{k=0}^n \sum_{j=1}^s W_{n-k,ij} Y_{kj},   \label{eq:rk16} \\
Y_{ni} = y_n + \tau \sum_{j=0}^s \alpha_{ij} Y'_{nj},  \qquad i=1,\ldots,s. \label{eq:rk17}
\end{align}
\end{lemma}
\noindent
The $y$ component of the Runge-Kutta discretization \eqref{eq:rk1}--\eqref{eq:rk4} of the system \eqref{eq:sys1}--\eqref{eq:sys2} thus again 
coincides with the approximation for the integro-differential equation \eqref{eq:idae}
obtained by the Runge-Kutta convolution quadrature scheme \eqref{eq:rk15}--\eqref{eq:rk17}.

\begin{remark}
For the Radau-IIA method with $s=1$ stages, one has $\A=\alpha_{11}=1$ and $\beta=\beta_1=1$. 
The symbol for the Runge-Kutta method then reads $\Delta(\xi)=(\frac{\xi}{1-\xi} + 1)^{-1} = 1-\xi$
which coincides with the associated polynomial $\delta(\xi)$ for the one-step BDF formula. 
In addition, also the weights $W_n=\omega_n$ coincide.
This can of course be expected, since both methods reduce to the implicit Euler method in that case. 
\end{remark}

\section{Details on the implementation} \label{sec:implementation}

For convenience of the reader, we now briefly discuss some aspects of the practical realization. We start by describing the efficient computation of the weights $\omega_n$ and $W_n$, and then briefly address the efficient evaluation of the convolution quadrature sums. 

\subsection{Computation of the quadrature weights for BDF formulas} 
To keep the presentation simple, we assume in the following that the linear subproblem \eqref{eq:sys2} is a single-input single-output system, i.e., matrices $B$ and $C$ have only one column. For multiple in- or outputs, the computations can simply be repeated for every single column. 
By choosing $\xi = \rho e^{i \phi}$ in equation \eqref{eq:bdf4},
we obtain 
\begin{align} \label{eq:weight1}
\sum_{n=0}^\infty \omega_n \rho^n e^{i n \phi} = \widehat K\left(\frac{\delta(\rho e^{i \phi})}{\tau}\right). 
\end{align}
Thus $\omega_n \rho^n$ are just the Fourier coefficients of the function $k(\phi) = \widehat K(\frac{\delta(\rho e^{i \phi})}{\tau})$. This immediately leads to the following explicit formula for the weights
\begin{align} \label{eq:weight2}
\omega_n = \frac{1}{2 \pi \rho^n} \int_0^{2\pi}   \widehat K\left(\frac{\delta(\rho e^{i \phi})}{\tau}\right) e^{-i n \phi} d\phi.
\end{align}
Computable approximations $\widetilde  \omega_n$ for the weights $\omega_n$ are then
obtained by replacing the integral with a numerical quadrature formula. 
Using the trapezoidal rule, one gets
\begin{align} \label{eq:trpz_weights}
\widetilde \omega_n = \frac{1}{L \rho^n}\sum_{l=0}^{L-1} 
\widehat K\left(\frac{\delta(\rho e^{i \phi_l})}{\tau}\right) 
e^{-i n \phi_l}, \qquad \phi_l=2 \pi l/L.
\end{align}
Under the assumption that $\widehat K(z)$ is analytic in a neighbourhood of the complex contour $z=\frac{\delta(\rho e^{i \phi})}{\tau}$, the quadrature error can be estimated by $|\omega_n - \widetilde \omega_n| \le a e^{-b L}$ for some $a,b>0$; see \cite{Henrici79} for details.
Let us note that a good approximation for the weights can therefore already be obtained for moderate $L$, 
i.e., only a reasonable number of evaluations of the transfer function will be required for setting up the method. 
\begin{remark} \label{rem:weights}
Following the discussion in \cite{Lubich88b}, an error bound 
$|\omega_n - \widetilde \omega_n| = O(\epsilon)$ can be obtained for the choice $\log \rho = O(h)$ and $L=O(\log \epsilon) N$ under rather general assumptions on the transfer function $\widehat K(z)$. 
Setting $L=O(N)$ and $\rho=\epsilon^{\frac{1}{2N}}$ still leads to a bound $|\omega_n - \widetilde \omega_n| = O(\sqrt{\epsilon})$. These choices will also be used in our numerical tests.
\end{remark}

\subsection{Computation of the weight matrices for Runge-Kutta formulas}

We again assume that $B$ and $C$ only consist of one single column.
By similar reasoning as in the previous section, we then obtain approximations 
\begin{align} \label{eq:rk-weights}
\widetilde W_n =  \frac{1}{L \rho^n}\sum_{l=0}^{L-1} 
\widehat K\left(\frac{\Delta(\rho e^{i \phi_l})}{\tau}\right) 
e^{-i n \phi_l}, \qquad \phi_l=2 \pi l/L,
\end{align}
for the weight matrices $W_n$ of the Runge-Kutta convolution quadrature method \eqref{eq:rk12}.

Before we proceed, let us also briefly comment on the computation of $\widehat K(M)$ for a matrix valued argument $M$.
Similar as above, we assume that the scalar valued transfer function $\widehat K(\xi) = \sum_{n=0}^\infty K_n \xi^n$ can be expanded into a power series, 
and we define 
\begin{align}
\widehat K(M) = \sum_{n=0}^\infty K_n M^n.
\end{align}
Now let $M=V \Lambda V^{-1}$ denote an eigenvalue decomposition of $M$ with diagonal matrix $\Lambda=\diag(\lambda_1,\ldots,\lambda_s)$ containing the eigenvalues. Then
\begin{align}
\widehat K(M) = V \widehat K(\Lambda) V^{-1} 
\qquad \text{with} \qquad \widehat K(\Lambda)=\diag(\widehat K(\lambda_1),\ldots,\widehat K(\lambda_s)).
\end{align}
Hence the evaluation of $\widehat K(\frac{\Delta(\xi)}{\tau})$ with matrix valued argument $\Delta(\xi) \in \CC^{s \times s}$ and $\xi \in \CC$ can be reduced to $s$ evaluations of $\widehat K(\lambda_k(\xi))$, $1 \le k \le s$ for a scalar valued argument and some elementary computations. 

\subsection{Remarks about the computational complexity} \label{sec:fast-summation}

A simple evaluation of the convolution sums in \eqref{eq:bdf9} or \eqref{eq:rk16} for all time steps $1 \le n \le N$ requires $O(N^2)$ operations. By an appropriate splitting of the sums and fast Fourier transforms, the complexity can however be reduced to $O(N \log N)$ operations; see \cite{Lubich88b,KapurEtAl96} for details. 

The online simulation of the coupled system \eqref{eq:sys1}--\eqref{eq:sys2} via the reduced problem \eqref{eq:idae} can therefore be obtained in $O(N \log N)$ complexity, if the size of the nonlinear subsystem \eqref{eq:sys1} is uniformly bounded. Up to the logarithmic factor, this is essentially the same computational cost as simulating the nonlinear subsystem \eqref{eq:sys1} alone. 
Using non-trivial modifications, the overall complexity could even be reduced to $O(N)$; see \cite{SchaedleEtAl06}.

\section{Field-circuit coupling} \label{sec:fieldcircuit}

As a typical application where our methods may be useful, we consider the simulation of electric circuits consisting of simple devices, which are described by simple algebraic relations, 
and complicated but linear {\em electromagnetic} components, which have to be modeled by Maxwell's equations.  
In this section, we present the basic model equations, and show that they perfectly fit into 
the class of problems considered in this paper. 
Two particular examples will be discussed in more detail in the following section. 

\subsection{Electric network model}

We assume that the network consists of an interconnection of resistors (R), capacitors (C), inductors (L), voltage sources (V), and an \emph{inductance like} electromagnetic element (M). 
Apart from the latter, all circuit elements shall be described by simple device relations \cite{GuentherEtAl05}. 
As usual, the interconnection of the different circuit elements is described by a directed graph whose (reduced) incidence matrix $A=[A_C,A_R,A_L,A_V,A_M]$ is naturally divided into subblocks for the different components. 
Using modified nodal analysis \cite{BartelBaumannsSchoeps11,GuentherEtAl05}, the time evolution of the circuit 
can then be described by the differential-algebraic equations
\begin{align}
A_C \ddt q_C(A_C^\top u) + A_R g_R(A_R^\top u) + A_L j_L + A_V j_V &= -A_M j_M  \label{eq:circuit1}\\
\ddt \phi(j_L) - A_L^\top u &= 0 \label{eq:circuit2}\\
A_V^\top u - v_s(t) &= 0.\label{eq:circuit3}
\end{align}
The node potentials $u$ and the currents $j_L$, $j_V$ through inductors and voltage sources are the unknowns of the system while the voltage source $v_s(t)$ describes the excitation. 
The vectors $v_X = A_X^\top u$ correspond to the voltage across the element $X$, and the nonlinear functions $q_C(v_C)$, $g_R(v_R)$, and $\phi_L(j_L)$ model the device characteristics of the capacitor, resistor, and inductor.
The current $j_M$ models the response of the electromagnetic element to the excitation by the voltage $v_M=A_M^\top u$; details will be given below.

Setting $y=(u,j_L,j_V)$, the system \eqref{eq:circuit1}--\eqref{eq:circuit3} attains the abstract form  \eqref{eq:sys1} with matrix
\begin{align*}
M(y)\!=\!\begin{pmatrix}A_C q_C'(A_C^\top u) A_C^\top & 0 & 0 \\ 0 & \phi_L'(j_L) & 0 \\ 0 & 0 & 0\end{pmatrix}
\end{align*}
and vectors 
\begin{align*}
f(y)\!=\!\begin{pmatrix} A_R g_R(A_R^\top u) + A_L j_L + A_V j_V \\                                                                                                                                                                                                                      -A_L j_L \\ A_V j_V - v_s(t)\end{pmatrix},
\qquad
C^\top z\!=\!\begin{pmatrix} -A_M j_M \\ 0 \\ 0 \end{pmatrix}.
\end{align*}
A complete definition of the vector $z$ will be given in the next section. 
The positive definite matrices $q_C'(v_C)$, $\phi_L'(j_L)$, and $g_R'(v_R)$ describe 
differential or incremental capacitance, inductance, and conductance of the corresponding devices.

We assume for now, and show in more detail below, that the constitutive equation for the electromagnetic element can be described in frequency domain by 
\begin{align} \label{eq:circuit4}
\widehat j_M(s) = \widehat k(s) \widehat v_M(s),  
\end{align}
where $\widehat k(s)$ is the reduced transfer function of the electromagnetic element. 
By inversion of the Laplace transform, the time domain response of the electromagnetic component can then be expressed as
\begin{align} \label{eq:circuit5}
j_M(t) = \int_0^t k(t-r) v_M(r) dr. 
\end{align}
The system \eqref{eq:circuit1}--\eqref{eq:circuit5} thus formally yields an integro-differential algebraic equation of the form \eqref{eq:idae}, as investigated in the previous sections. 
In the following, we discuss in more detail the derivation of the transfer function $\widehat k(s)$ for the electromagnetic element.

\subsection{Electromagnetic device model}

\begin{figure}[ht!]

\centering
\begin{tikzpicture}[even odd rule]
\begin{scope}[spy using outlines=
	{circle, magnification=3, size=4cm, connect spies}]
\draw (0,0) rectangle (4,4);
\draw [fill = gray!20] (2,2) circle[radius = 2/3] circle[radius = 4/3];
\draw [black] (2.47, 2.47) -- (2.94,2.94);
\node at (3.2, 3.2) {$\Sigma$};
\node at (0.6, 3.4) {$\Omega_{nc}$};
\node at (2.1, 1.) {$\Omega_c$};
\node at (0.6, 1.15) {$\Gamma$};
\spy [black] on (2.7,2.7)
   in node[fill=white] at (7,2);
\end{scope}
\node at (6.4, 2.5) {\huge $\Sigma_1$};
\node at (7.5, 1.2) {\huge $\Sigma_0$};
\draw [dashed, black, thick] (6.25,1.4) -- (7.65,2.8);
\draw [dashed, black, thick] (6.4,1.25) -- (7.8,2.65);
\draw [thick] (7.475, 2.625) to[out=140,in=180]  (8.5,4) to[out=0,in=180]  (10,3.5);
\draw [thick] (7.325, 2.175) to[out=-40,in=180] (10,0.5);
\draw (10,0.5) to[vco,o-o,l_= $v_M(t)$] (10,3.5);
\end{tikzpicture}
\caption{Schematic view of the domain for the electromagnetic element.}
\label{fig:schematic}
\end{figure}
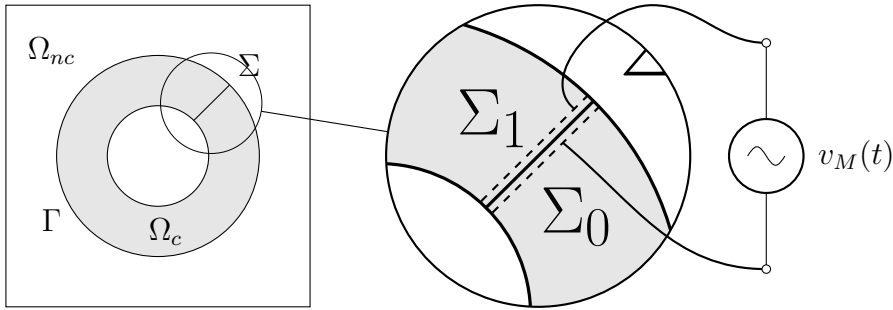

For ease of presentation, we assume that the electromagnetic element (M) represents a solid conductor $\Omega_c$  
embedded in a non-conducting region $\Omega_{nc}=\Omega\setminus\Omega_c$. 
Stranded conductors will also be discussed in the numerical tests below. 
The interface between the two regions is denoted by $\Gamma=\overline\Omega_c \cap \overline\Omega_{nc}$. 
A voltage $v_M$ is applied across the interface $\Sigma$, whose two sides $\Sigma_0$, $\Sigma_1$ model two electric contacts via which the electromagnetic element is connected to the electric circuit; see Figure~\ref{fig:schematic} for a sketch of the geometry.

Under the assumption that the magneto-quasistatic approximation is valid, 
the evolution of the  magnetic field can then be described by \cite{AlonsoValli10,HiptmairSterz05}
\begin{align} \label{eq:eddy1}
 \int_{\Omega_c} \sigma \ddt a(t) \cdot a'(t) dx + \int_\Omega \nu \curl a(t) \cdot \curl a' dx = -v_M(t) \int_{\Omega_c} \sigma p \cdot a' dx.
\end{align}
Here $a(t)$ denotes the magnetic vector potential which is element of an apprioriate subspace $V$ of $H(\curl;\Omega)$ incorporating gauging and boundary conditions, $\sigma$ is the electric conductivity, and $\nu = \mu^{-1}$ the inverse of magnetic permeability. The above variational equation is assumed to hold for all test function $a' \in V$ and all $t>0$,
and we further assume that $a(0) \equiv 0$ at time $t=0$ in the sequel.
The function $p$ in \eqref{eq:eddy1} can be chosen as $p=\nabla \phi$, where $\phi \in H^1(\Omega_c \setminus \Sigma)$ denotes a scalar potential with $\phi|_{\Sigma_0}=0$ and $\phi|_{\Sigma_1}=1$. 
The total current through the contacts $\Sigma_0$, $\Sigma_1$ arising in response to the excitation by the voltage $v_M(t)$ can finally be expressed as
\begin{align} \label{eq:eddy2}
j_M(t) = \int_{\Omega_c} \sigma \ddt a(t) \cdot p \; dx + v_M(t) \int_{\Omega_c} \sigma |p|^2 dx.  
\end{align}
We refer to \cite[Sec.~5.2]{HiptmairSterz05} for details on this particular model and to \cite{Dular04,HiptmairSterz05} for alternative formulations and other geometric configurations.

A semi-discretization of \eqref{eq:eddy1}--\eqref{eq:eddy2} in space by appropriate $H(\curl)$-conforming finite elements leads to a differential-algebraic system of the form 
\begin{align}
M_\sigma \ddt a(t) + K_\nu a(t) &= -B_1 v_M(t) \label{eq:eddy4}\\
B_1^\top \ddt a(t) - j_M(t)     &= -B_2 v_M(t). \label{eq:eddy5}
\end{align}
The current $j_M$ resulting in response to the excitation of $v_M$ can then be expressed in frequency domain by $\widehat j_M(s) = \widehat k(s) \widehat v_M(s)$ with
\begin{equation}
\label{eq:transfer_reduced}
\widehat k(s) = B_2 - s B_1^\top (s M_\sigma  + K_\nu)^{-1} B_1.
\end{equation}

Let us finally illustrate that the above system falls into the class of problems discussed in the paper. 
Recall that $y=(u,j_L,j_V)$ and thus $v_M = A_M^\top u = \begin{pmatrix} A_M^\top & 0 & 0 \end{pmatrix} y$.  
Then by setting $z=(a,j_M)$, the system \eqref{eq:eddy4}--\eqref{eq:eddy5} can be cast into the abstract form \eqref{eq:sys2} with 
\begin{align*} 
E=\begin{pmatrix} M_\sigma & 0 \\ B_1^\top & 0 \end{pmatrix}, \quad 
A=\begin{pmatrix} K_\nu & 0 \\ 0 & -1\end{pmatrix}, \quad 
\text{and} \quad 
B=\begin{pmatrix} -B_1 \\ -B_2 \end{pmatrix} \begin{pmatrix} A_M^\top & 0 & 0 \end{pmatrix},
\end{align*}
and the response of this linear subsystem can be expressed as
\begin{align} \label{eq:transferK} 
C^\top \widehat z(s) 
= \widehat K(s) \widehat y(s) 
= - \begin{pmatrix} A_M \\ 0 \\ 0 \end{pmatrix} \widehat k(s) \begin{pmatrix} A_M^\top & 0 & 0 \end{pmatrix} \widehat y(s).  
\end{align}
Let us note that only the transfer function $\widehat k(s)$ and the matrices $A_M$, which can be accessed efficiently, will be required for the further computations.

\subsection{Coupled electric circuit--electromagnetic device model}

As shown above, the coupled field-circuit problem \eqref{eq:circuit1}--\eqref{eq:circuit3} and \eqref{eq:eddy4}--\eqref{eq:eddy5} perfectly fits into the structure of problems investigated in the previous sections. 
Further note that under general assumptions on network topology and the constitutive equations of the individual devices, the coupled differential-algebraic system can be shown to have index-1; see e.g. \cite{CortezEtAl18,BartelBaumannsSchoeps11}. One can thus expect full convergence orders for the BDF and Radau-IIA methods discussed in the previous sections; we refer to \cite{BrenanCampbellPetzold96,HairerWanner10b} for details. 

\section{Numerical tests} \label{sec:numerics}

In the following, we illustrate our theoretical results by some numerical tests. 
We start by considering a simple linear time invariant model problem, and then briefly discuss the application to a linear field--nonlinear circuit model describing a half rectifier \cite{SchoepsEtAl10}.

\subsection{Model problem}

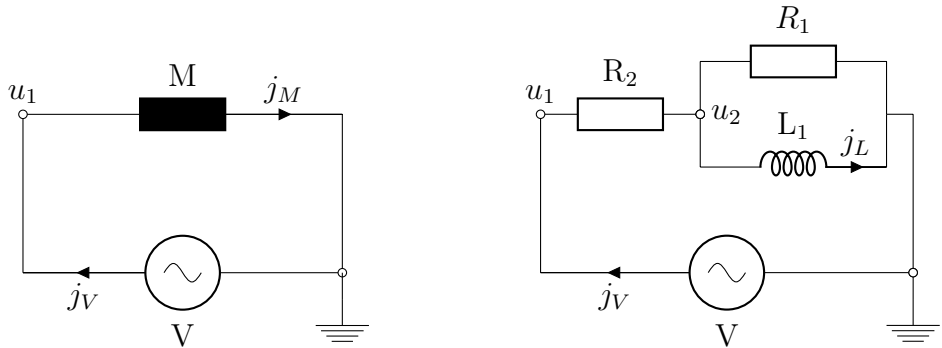
\begin{figure}[ht!]
\centering
\begin{subfigure}[b]{0.45\linewidth}
\centering
\begin{circuitikz}[scale = 0.7]
\draw
(0,3) to[short,-o] (0,6) node[above]{$u_1$} 
to [fullgeneric=M, i=$j_M$] (6,6) to [short,-] (6,3) 
to [vco,o-,l^=V,i=$j_V$] (0,3);
\draw
(6,3) to[short, -] (6,2.5) node[ground]{};
\end{circuitikz}
\end{subfigure}
\vspace*{2em}
\begin{subfigure}[b]{0.45\textwidth}
\centering
\begin{circuitikz}[scale = 0.7]
\draw (0,-0.5) -- (0,2.5) node[above]{$u_1$} to [R=R$_2$,o-o](3,2.5) node[right]{$u_2$};
\draw (3,2.5) to[short,o-] (3,3.5) to [R=$R_1$,-,] (6.5,3.5) to[short,-] (6.5,2.5) ;
\draw (3,2.5) to[short,o-] (3,1.5) to [L=L$_1$,-,i=$j_L$] (6.5,1.5) to[short,-] (6.5,2.5);
\draw (6.5,2.5) to [short, -] (7,2.5) -- (7,-0.5);
\draw (7,-0.5) to [vco,l^=V,o-,i=$j_V$] (0,-0.5);
\draw (7,-0.5) to[short, o-] (7,-1) node[ground]{};
\end{circuitikz}
\end{subfigure}
\caption{\label{fig:modelproblem}Electric network for the model problem with an electromegnetic element (left) and an equivalent circuit (right).}
\end{figure}

Let us start by considering a simple circuit consisting only of a voltage source and one electromagnetic element. For comparison with other model reduction approaches, we will also briefly discuss the replacment of the electromagnetic element by an equivalent circuit; see Figure~\ref{fig:modelproblem} for a sketch of the network topology. 
\subsubsection{Electromagnetic device}
As a model for the electromagnetic element, we consider a two-dimensional setup in which the vector potential is of the form $a = (a_x,a_y,0)$ with components $a_x,a_y$ independent of the $z$-coordinate. We choose $\Omega = (-1,1)^2$ as computational domain with conductor loop $\Omega_c$ between inner and outer radius $r_{in} = 1/3$ and $r_{out} = 2/3$; see Figure~\ref{fig:schematic} for an illustration. The permeability is set to $\mu \equiv 1$ in the whole domain $\Omega$, and we choose $\sigma = 1$ in $\Omega_c$ and $\sigma = 0$ in $\Omega_{nc}$ for the conductivity in the subdomains. 
The voltage $v_M$ is applied across the interface $\Sigma$, resulting in a current flowing in circles through the conducting domain $\Omega_c$; cf. Figure~\ref{fig:schematic}. 
As can be deduced from equation \eqref{eq:eddy2}, this current will additionally be 
altered by induction effects.

For the discretization of the electromagnetic field equations \eqref{eq:eddy1}--\eqref{eq:eddy2}, we 
utilize second order Nedelec finite elements on a quadrilateral mesh with curved elements. The system matrices for \eqref{eq:eddy4}--\eqref{eq:eddy5} were assembled in the C++ library \texttt{Concepts} \cite{Concepts}.

\subsubsection{Convergence plots}

\begin{figure}[ht!]
\centering
\begin{tikzpicture}[scale=0.70]
\begin{loglogaxis}[
width = 0.6\linewidth, 
height = 0.5\linewidth,
legend style={at={(1,0)},anchor=south east},
grid=both,
minor grid style={gray!25}, 
major grid style={gray!25},
width=\linewidth]
\addplot[orange, only marks, mark=o,mark options={scale=2}] table[x={dt}, y={CQ}] {conv_BDF_1.dat};
\addlegendentry{BDF-1};
\addplot[purple, only marks, mark=o,mark options={scale=2}] table[x={dt}, y={CQ}] {conv_BDF_2.dat};
\addlegendentry{BDF-2};
\addplot[red, only marks, mark=o,mark options={scale=2}] table[x={dt}, y={CQ}] {conv_RK_2.dat};
\addlegendentry{RadauIIA-2};
\addplot[blue, only marks, mark=o,mark options={scale=2}] table[x={dt}, y={CQ}] {conv_RK_3.dat};
\addlegendentry{RadauIIA-3};
\addplot[orange, thick, only marks, mark=x, mark options={scale=2}] table[x={dt}, y={RK}] {conv_BDF_1.dat};
\addplot[purple, thick, only marks, mark=x, mark options={scale=2}] table[x={dt}, y={RK}] {conv_BDF_2.dat};
\addplot[red, thick, only marks, mark=x, mark options={scale=2}] table[x={dt}, y={RK}] {conv_RK_2.dat};
\addplot[blue, thick, only marks, mark=x, mark options={scale=2}] table[x={dt}, y={RK}] {conv_RK_3.dat};
\addplot[dashed,domain= 0.002:0.5]{0.02*x^1};
\addplot[dashed,domain= 0.002:0.5]{0.012*x^2};
\addplot[dashed,domain= 0.002:0.5]{0.012*x^3};
\addplot[dashed,domain= 0.002:0.5]{0.005*x^5};
\end{loglogaxis}
\end{tikzpicture}
\caption{Errors in the $y$ component vs. time step $\tau$ for the one- and two-step BDF and the Radau-IIA methods with two and three stages applied to the coupled system (crosses) and the reduced problem (circles). The dashed lines are the theoretical convergence rates.}
\label{fig:convergence}
\end{figure}
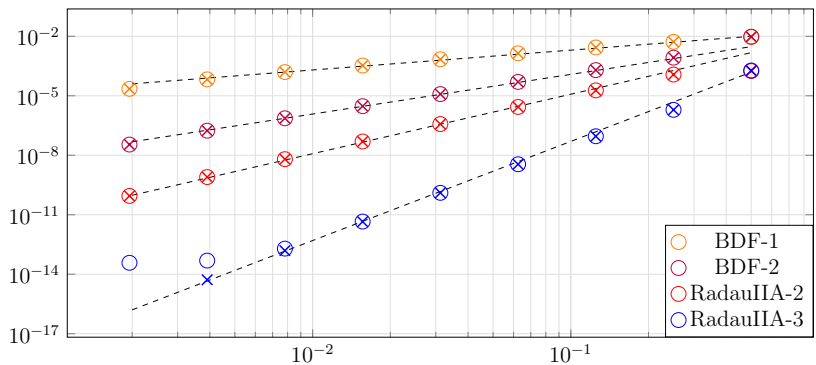

To illustrate the theoretical results of paper, we first report on results for an excitation of the model problem with source current $v_s(t)=\sin(3/2\pi t)$ on a time interval $t \in [0,1]$. 
In Figure~\ref{fig:convergence}, we display the simulation errors 
for BDF and Radau-IIA methods of various orders applied to the coupled system \eqref{eq:sys1}--\eqref{eq:sys2} and to the corresponding methods with convolution quadrature applied to the integro-differential equation \eqref{eq:idae}. 
A fixed spatial discretization with 19979 degrees of freedom was used for magnetic field computation.
Following Remark~\ref{rem:weights}, we chose $L=3 N$ and $\rho = \exp(-\tau)$ in the formulas \eqref{eq:trpz_weights} and \eqref{eq:rk-weights} for the computation of the quadrature weights. 
Hence we can expect to obtain exact weights up to machine precision.

As predicted by our theoretical results, the convolution quadrature methods applied to the reduced problems \eqref{eq:idae} yield the same convergence rates as the discretizations of the coupled problem \eqref{eq:sys1}--\eqref{eq:sys2} by standard one- and multistep methods. In fact, the difference between the corresponding numerical solutions of the two equivalent formulations was in the order of round-off errors in all tests, which can only be seen in the results for the three-step Radau-IIA method at very small step sizes. 

Let us emphasize that Figure~\ref{fig:convergence} illustrates only the time discretization errors, i.e., 
the computations were done for a fixed space discretization. By comparison with computations on 
a refined mesh, the spatial discretization error was estimated to be $O(10^{-8})$.

\subsubsection{Transfer function and approximation by an equivalent circuit}

\begin{figure}[ht]
\centering
\medskip
\begin{subfigure}{0.45\linewidth}
\centering
\begin{tikzpicture}
\footnotesize
\begin{semilogxaxis}[
xlabel={$\omega$ (rad/s)},
ylabel={$| \widehat{k}(i\omega) |$ (dB)},
grid=both, 
minor grid style={gray!25}, 
major grid style={gray!25},
no marks,
width=\linewidth]
\addplot[red, thick] table[x={w}, y={mag}] {trans.dat};
\addplot[black, dashed] table[x={w}, y={mag}] {trans_approx.dat};
\end{semilogxaxis}
\end{tikzpicture}
\end{subfigure}%
\hfill
\begin{subfigure}{0.45\linewidth}
\centering
\begin{tikzpicture}
\footnotesize
\begin{semilogxaxis}[
xlabel= $\omega$ (rad/s),
ylabel= $\arg \widehat{k}(i\omega) $ (deg) ,
grid=both, 
minor grid style={gray!25}, 
major grid style={gray!25},
no marks,
width=\linewidth]
\addplot[red, thick] table[x={w}, y={phase}] {trans.dat};
\addplot[black, dashed] table[x={w}, y={phase}] {trans_approx.dat};
\end{semilogxaxis}
\end{tikzpicture}
\end{subfigure}%
\caption{Bode diagram for the transfer functions of the electromagnetic 
element (red line) and the equivalent circuit (black dashed).}
\label{fig:bodeplot}
\end{figure}
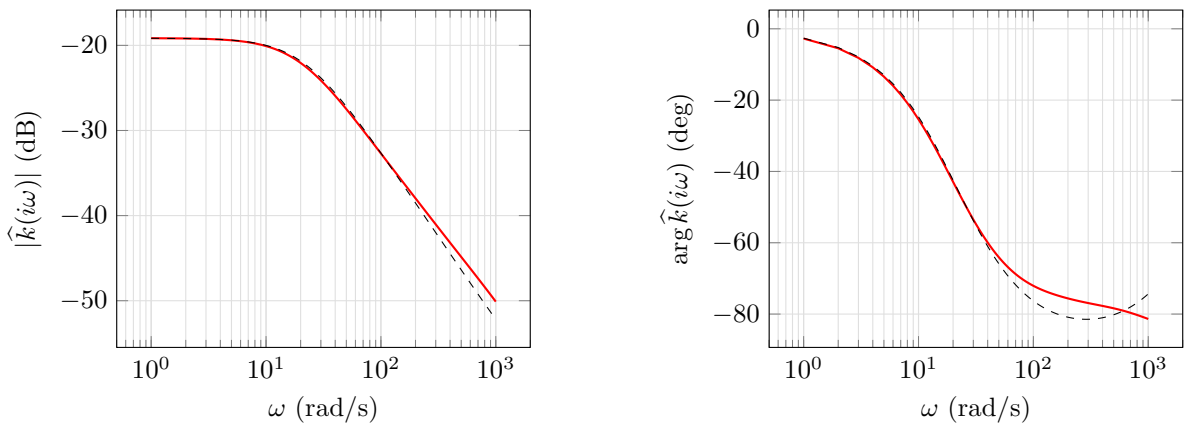

For comparison with other model reduction approaches, we now briefly discuss the approximation of the electromagnetic element by an equivalent circuit; see Figure~\ref{fig:modelproblem}.
Let us start by recalling that the response of the electromagnetic field element can be described
in frequency domain by $j_M(s) = \widehat k(s) \widehat v_M(s)$ with 
\begin{align} \label{eq:response}
\widehat k(s) =   B_2 - s B_1^\top (s M_\sigma + K_\nu)^{-1} B_1.
\end{align}
Mimicking the algebraic form of this formula, we consider a simple (1,1)-rational approximation of the transfer function \eqref{eq:response} given by 
\begin{align} \label{eq:responseEC}
\widehat k_{EC}(s) = a - \frac{s}{c s + d}.
\end{align}
Let us note that $\widehat k_{EC}(s)$ just amounts to the transfer function of the equivalent circuit depicted in Figure~\ref{fig:modelproblem} (right), with parameters $R_2 = 1/a$, $R_1 = 1/(ca^2-a)$ and $L_1 = 1/(a^2d)$. 
For our test, we determined values $R_1=1623$, $R_2=9.08$, and $L_1=0.4174$ by fitting $\widehat k_{EC}(s)$ to $\widehat k(s)$ via the AAA algorithm \cite{AAA18}. 
Thousand evaluations of the transfer function were used to obtain a reliable approximation. 
In Figure~\ref{fig:bodeplot}, we compare the exact transfer function with the corresponding 
rational fit in a Bode diagram. 

One can deduce that, at least for sufficiently small frequencies, the response of the simple electromagnetic device can be approximated well in amplitude and in phase by that of an appropriate equivalent circuit.
Let us note, however, that even for a low frequency input signal $v(t) = \sin(3/2 \pi t)$, used in our simulations above, a systematic error of $O(10^{-4})$ was observed in the corresponding time-domain simulations. 
To improve the accuracy, higher order rational approximation have to be used
for which the interpretation as equivalent circuit becomes more difficult. 

\subsubsection{Computational performance}

We next discuss the computational complexity of the convolution quadrature method
and compare it to more standard approaches.

\subsubsection*{Offline phase}

The computation of the weights via \eqref{eq:trpz_weights} or \eqref{eq:rk-weights} 
is the compute intensive part of our approach, but this can be done in an offline stage. 
Recall that $L = O(N)$ was chosen in our computations, i.e., the number of evaluations of the transfer function is comparable to the number of time steps. Since every evaluation involves the solution of a linear system for the electromagnetic field part, the complexity is comparable to that for one solution of the system \eqref{eq:sys1}--\eqref{eq:sys2} with an implicit Runge-Kutta or multistep method not explicitly 
making use of the time-invariant nature of the problem. 
Let us further note that the complexity of computations in the offline phase could be reduced substantially by using \textit{fast and oblivious convolution quadrature} \cite{SchaedleEtAl06}, which only requires about $O(\log N)$ evaluations of the transfer function to set up the quadrature weights.

Using the time invariance of the model problem under consideration and assuming a fixed time step $\tau$, 
the system matrices for the implicit Runge-Kutta and multistep methods applied to \eqref{eq:sys1}--\eqref{eq:sys2} 
will not change during simulation. They can thus be factorized once prior to the 
computation. 

\subsubsection*{Online phase}

After the system analysis performed in the offline stage has been performed, 
we now compare the computation times required for the actual simulation of the coupled problem \eqref{eq:sys1}--\eqref{eq:sys2} and the reduced model \eqref{eq:idae}. As before, we use a fixed spatial discretization with 19970 degrees of freedom for the finite element approximation of the electromagnetic element. 
The results of our computations are given in Table~\ref{tab:complexity}. 
\begin{table}[htp] \small
  \centering
  \setlength\tabcolsep{1em}
  \begin{tabular}{r|r|r|r}
	$N$  & RK-coupled     & RK-CQ-red  & RK-equiv \\
    \hline
     16  &    1.265 & 0.008  & 0.007 \\
     32  &    2.452 & 0.015  & 0.012 \\
     64  &    4.933 & 0.032  & 0.024 \\
    128  &    9.732 & 0.095  & 0.029 \\
    256  &   19.365 & 0.318  & 0.057 \\
    512  &   39.280 & 1.185  & 0.115 \\
\end{tabular}
\medskip
\caption{Online computation times (sec) for three-stage Radau-IIA method applied to the coupled system (RK-coupled), the reduced model (RK-CQ-red), and the equivalent circuit model (RK-equiv); cf. Figure~\ref{fig:schematic}.}
\label{tab:complexity}
\end{table}

Let us note that we used a simple implementation of the convolution quadrature sums here and, therefore,
the complexity of the Runge-Kutta convolution quadrature method is $O(N^2)$, while that for the 
other two methods is $O(N)$. 
Following our remarks in Section~\ref{sec:fast-summation}, 
the online complexity of our algorithm could be reduced to $O(N \log N)$.  
Therefore, only the computation times for small $N$ allow for a fair comparison here. 

The online simulation times for the reduced problem \eqref{eq:idae} by the proposed Runge-Kutta convolution quadrature method thus is essentially the same as that for the equivalent circuit model. Recall, however, that the systematic error in the latter was about $10^{-4}$, 
while no systematic error is present in the proposed Runge-Kutta convolution quadrature method, except 
round-off and spatial discretization errors. 

Let us further note that we explicitly utilized the time invariance of the problem also in the implementation of the Runge-Kutta method for the coupled field problem here, i.e., the system matrix was factorized once before simulation, and only forward-backward substitutions were performed at every 
time step in the online phase.  
Nevertheless, even for the simple two dimensional electromagnetic field problem, 
the simulation times for the coupled field-circuit problem were significantly higher than those for the reduced models. This difference will even become bigger when larger pde models are required to approximate the electromagnetic element with sufficient accuracy.

In summary, we can thus say that the Runge-Kutta convolution quadrature method allows an efficient online 
simulation of the coupled electromagnetic field-circuit problem, essentially at the cost of simulating 
the circuit model only, while delivering accuracy and stability of the simulation of the 
fully coupled model.

\subsection{Simulation of a half-rectifier}

\begin{figure}[ht!]
\centering
\begin{subfigure}[t]{0.70\linewidth}
\centering
\pgfdeclarelayer{bg}    
\pgfsetlayers{bg,main}  
\begin{circuitikz}[scale = 0.6]
\draw (0,0) node[transformer core, yscale=1.5] (T) {}
(T.A1) node[anchor=east] {}
(T.A2) node[anchor=east] {}
(T.B1) node[anchor=west] {}
(T.B2) node[anchor=west] {}
(T.base) node{};
\begin{pgfonlayer}{bg}  
\draw [dashed,fill=gray!10] ($(T.A2)-(0.3,0.3)$)  node[below right]{$M$} rectangle ($(T.B1)+(0.3,0.3)$);
\end{pgfonlayer}
\draw (T.A2) -- ($(T.A2)-(2,0)$) to[vco,l^=$V$,o-o] ($(T.A1)-(2,0)$) node[above]{$u_1$} -- (T.A1); 
\draw (T.B1) -- ($(T.B1)+(2,0)$) node[above]{$u_2$} to[C,l=$C$,o-] ($(T.B2)+(2,0)$) -- (T.B2);
\draw ($(T.B1)+(2,0)$) to[Do,l=$D$,-o] ($(T.B1)+(6,0)$) node[above]{$u_3$} to[R,l=$R$] ($(T.B2)+(6,0)$) to[short, -] ($(T.B2)+(2,0)$);
\draw ($(T.A2)+(1.05,0)$) to[short,-] ($(T.B2)-(1.05,0)$);
\draw ($(T.A2)-(2,0)$) to[short,-]  ($(T.A2)-(2,0.5)$) node[ground]{};
\end{circuitikz}
\end{subfigure}
\begin{subfigure}[t]{0.25\linewidth}
\begin{tikzpicture}[scale = 0.072]
\draw (-20,-25) -- (20,-25) -- (20,25) -- (-20,25) -- (-20,-25);
\draw [draw=black, fill=gray, opacity =0.2] (-12,-20) -- (12,-20) -- (12,20) -- (-12,20) -- (-12,-20);
\draw (-12,-20) -- (12,-20) -- (12,20) -- (-12,20) -- (-12,-20);
\draw [draw=black, fill=white](-4,-12) -- (4,-12) -- (4,12) -- (-4,12) -- (-4,-12);
\draw (-12,-11) -- (-15,-11) -- (-15,11) -- (-12,11);
\draw (12,-11) -- (15,-11) -- (15,11) -- (12,11);
\draw (-4,-11) -- (-1,-11) -- (-1,11) -- (-4,11);
\draw (4,-11) -- (1,-11) -- (1,11) -- (4,11);
\node [align=center,] at (-13.5,0) {\tiny$\otimes$};
\node [align=center,] at (-13.5,6) {\tiny$\otimes$};
\node [align=center,] at (-13.5,-6) {\tiny$\otimes$};
\node [align=center,] at (2.5,0) {\tiny$\otimes$};
\node [align=center,] at (2.5,6) {\tiny$\otimes$};
\node [align=center,] at (2.5,-6) {\tiny$\otimes$};

\node [align=center,] at (13.5,0) {\tiny$\odot$};
\node [align=center,] at (13.5,6) {\tiny$\odot$};
\node [align=center,] at (13.5,-6) {\tiny$\odot$};
\node [align=center,] at (-2.5,0) {\tiny$\odot$};
\node [align=center,] at (-2.5,6) {\tiny$\odot$};
\node [align=center,] at (-2.5,-6) {\tiny$\odot$};
\draw (13.5,11) |- (26,22);
\draw (13.5,-11) |- (26,-22);
\draw (-13.5,11) |- (-26,22);
\draw (-13.5,-11) |- (-26,-22);
\node at (0,-32) {};
\end{tikzpicture}
\end{subfigure}
\caption{Sketch of the rectifier circuit (left) and geometry of the transformer (right) modeled by the electromagnetic field equations.}
\label{fig:rectifier}
\end{figure}
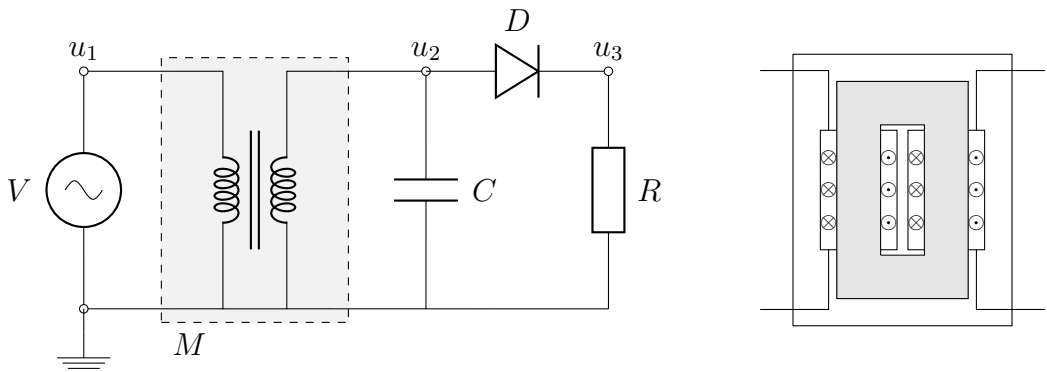

As a second and more practical example, we choose a test problem discussed in \cite{SchoepsEtAl10}
which is concerned with the simulation of a half-rectifier. A schematic sketch of the corresponding 
circuit is depicted in Figure~\ref{fig:rectifier}.

The electromagnetic component is represented by a 2D model of the transformer depicted in Figure~\ref{fig:rectifier} under the assumption that the magnetic vector potential is of the form $a=(0,0,a_z)$ with $a_z=a_z(x,y)$; the current is injected to the transform through stranded conductors. 
We refer to \cite[Sec.~3.2]{HiptmairSterz05} for derivation of the model equations, and to \cite[Sec.~6.2]{Schoeps11Diss} for material properties and details on transformer model. 
Discretization by finite elements then lead to a differential algebraic system of the form
\begin{align}
M_\sigma \ddt a(t) + K_\nu a(t) - B_3^\top j_M(t) &= 0 \label{eq:eddy8}\\
B_3 \ddt a(t)  &= v_M(t). \label{eq:eddy9}
\end{align}
Let us note that the transformer represents a system with two input and two output ports.
The transfer function $\widehat{k}(s)$ at frequency $s$ thus is a $2 \times 2$ matrix 
defined by   
\begin{align} \label{eq:eddy10}
\widehat j_M(s) 
= \widehat k(s) \widehat v_m(s)
= (s B_3 (s M_\sigma + K_\nu)^{-1} B_3)^{-1} \widehat v_M(s).
\end{align}
We again use \eqref{eq:transferK} to define the transfer function of the linear subsystem 
given in \eqref{eq:Ctzh}.

All circuit elements, except the diode, are assumed to have a linear response with $C = 10^{-12}$ and $R =10000$. The nonlinear voltage-current relation for the diode is given by $j_D = 2.5\cdot 10^{-6} (\exp (4 v_D)+1)$, which amounts to that of a Shockley diode.

For our numerical tests, we apply $V(t) = 250\sin (5\pi t)$ as an input signal and consider time interval $t \in [0,1]$ with $N =1000$ points for time stepping. 
The convolution quadrature weights are computed by formula \eqref{eq:trpz_weights} with $L=N$ and $\rho = \sqrt[2N]{\epsilon}$, where $\epsilon = 10^{-16}$; see Remark \ref{rem:weights}.
In Figure~\ref{fig:rectifier}, we display the node potentials $u_1 (t)$ and $u_3(t)$, which correspond to the input and output voltages of the rectifier circuit. 
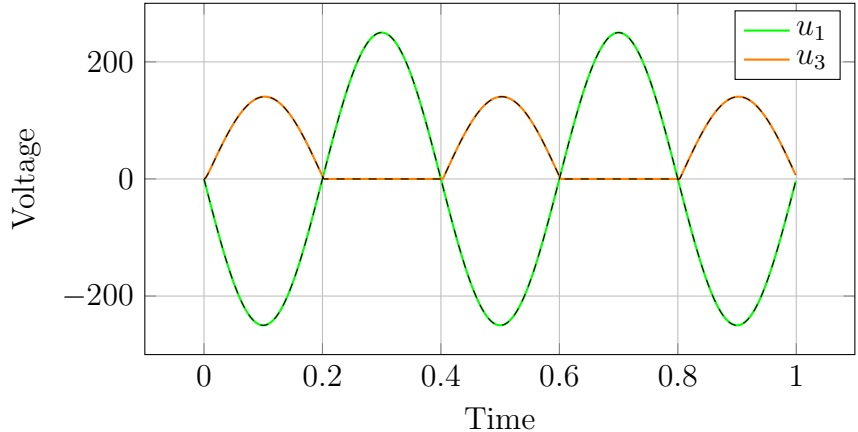
\begin{figure}[ht!]
\centering
\begin{tikzpicture}[scale=1]
\begin{axis}[
width = 0.7\linewidth, 
height = 0.4\linewidth,
ylabel={Voltage},
xlabel={Time}, 
grid=both]
\addplot[green, thick] table[x={dt}, y={CQ1}] {rectifier.dat};
\addlegendentry{$u_1$};
\addplot[orange, thick] table[x={dt}, y={CQ3}] {rectifier.dat};
\addlegendentry{$u_3$};
\addplot[dashed] table[x={dt}, y={RK1}] {rectifier.dat};
\addplot[dashed] table[x={dt}, y={RK3}] {rectifier.dat};
\end{axis}
\end{tikzpicture}
\caption{Node potentials $u_1$ and $u_3$ for coupled (dashed) and reduced (colored) models obtained by implicit Euler method.}
\label{fig:rectifier_result}
\end{figure} 

The simulation results show the expected behavior, i.e., the diode hinders the current to flow in the wrong direction, which results in a zero potential $u_3$ during half of the period. In the second half of 
the period, the diode lets the current pass freely and the output voltage $u_3$ amounts to minus the input voltage. This behavior also explains the name \emph{half-rectifier}.
Let us note that, in accordance with our theorems, the simulation of the coupled system \eqref{eq:sys1}--\eqref{eq:sys2} and the reduced system \eqref{eq:idae} yield identical results.
The latter, however, only requires the time integration of the nonlinear subsystem, which here is of dimension four, while the first approach needs to solve also for the linear subsystem \eqref{eq:sys2} 
in every time step. 

\section*{Acknowledgements}

The authors are grateful for financial support by the ``Excellence Initiative'' of the German Federal and State Governments via the Graduate School of Computational Engineering GSC~233 at Technische Universität Darmstadt and by the German Research Foundation (DFG) via grants TRR~146, TRR~154, and Eg-331/1-1.

\end{document}